\documentclass[aap,MSNbibl,nameyear,dvips]{arximspdf}

\doi{10.1214/11-AAP775}
\volume{21}
\issue{5}
\pubyear{2011}
\firstpage{2050}
\lastpage{2051}

\makeatletter

\newtheorem{cor}{Corollary}[section]

\makeatother

\begin{document}
\begin{frontmatter}

\title{Corrections\\
Occupation and local times for skew Brownian motion with applications to dispersion across an~interface}\vspace*{4pt}
\runtitle{Corrections}
\pdftitle{Corrections
Occupation and local times for skew Brownian motion with applications
to dispersion across an interface}

\textit{Ann. Appl. Probab.} \textbf{21} (2011) 183--214

\begin{aug}
\author[A]{\fnms{Thilanka} \snm{Appuhamillage}\corref{}\ead[label=e1]{ireshara@math.oregonstate.edu}},
\author[A]{\fnms{Vrushali} \snm{Bokil}\ead[label=e2]{bokilv@math.oregonstate.edu}},
\author[A]{\fnms{Enrique} \snm{Thomann}\ead[label=e3]{thomann@math.oregonstate.edu}},
\author[A]{\fnms{Edward} \snm{Waymire}\ead[label=e4]{waymire@math.oregonstate.edu}}
and
\author[B]{\fnms{Brian} \snm{Wood}\ead[label=e5]{brian.wood@oregonstate.edu}}
\runauthor{T. Appuhamillage et al.}
\affiliation{Oregon State University}
\address[A]{T. Appuhamillage\\
V. Bokil\\
E. Thomann\\
E. Waymire\\
Department of Mathematics\\
Oregon State University\\
Corvallis, Oregon 97331\\
USA\\
\printead{e1}\\
\phantom{E-mail: }\printead*{e2}\\
\phantom{E-mail: }\printead*{e3}\\
\phantom{E-mail: }\printead*{e4}}
\address[B]{B. Wood\\
School of Chemical, Biological \\
\quad and Environmental Engineering \\
Oregon State University\\
Corvallis, Oregon 97331\\
USA\\
\printead{e5}}
\end{aug}

\received{\smonth{3} \syear{2011}}


%
\begin{keyword}[class=AMS]
\kwd[Primary ]{60K35}
\kwd{60K35}
\kwd[; secondary ]{60K35}.
\end{keyword}
\begin{keyword}
\kwd{Skew Brownian motion}
\kwd{advection-diffusion}
\kwd{local time}
\kwd{occupation time}
\kwd{elastic skew Brownian motion}
\kwd{stochastic order}
\kwd{first passage time}.
\end{keyword}

\end{frontmatter}

\vspace*{-6pt}
The nonnegative parameter $\gamma= |(2\alpha-1)v|$ appearing in
the formula for the transition probabilities for $\alpha$-skew Brownian
motion with drfit $v$
in Theorem 1.3
should be replaced by the parameter $\gamma= (2\alpha-1)v\in(-v,v)$.
The formula is then correct, as can be checked by (tedious)
differentiations for the backward equation and
interface condition in the backward variable $x$.
However, it is only in the cases when $(2\alpha-1)v \ge0$ that the
probabilistic interpretation
of $\gamma$ as a skew-elasticity parameter applies.

In addition, the corrected
display to Corollary \ref{corc33} that follows from integration of the formula
in Corollary 1.2
giving the trivariate density is as follows:
\vspace*{-9pt}

\setcounter{cor}{2}
\setcounter{section}{3}
\begin{cor}\label{corc33}
If $x\geq0$ we have
\begin{eqnarray*}
&&
P_x\bigl(B_t^{(\alpha)}\in dy, \ell_t^{(\alpha)}\in d\ell\bigr)\\[-3pt]
&&\qquad=\cases{
\displaystyle \frac{2(1-\alpha)(l-y+x)}{\sqrt{2\pi t^3}}
\exp\biggl\{-\frac{(l-y+x)^2}{2t}\biggr\}\,dy\,dl,
\vspace*{-8pt}
\cr
\hspace*{115.8pt}\hspace*{70.6pt}\qquad\qquad\mbox{if } y\le0,l \ge0, \vspace*{5pt}\cr
\displaystyle \frac{2\alpha(l+y+x)}{\sqrt{2\pi t^3}}\exp\biggl\{-\frac{(l+y+x)^2}{2t}\biggr\}
\,dy\,dl\vspace*{1pt}
\cr
\qquad{} + \displaystyle \frac{1}{\sqrt{2\pi
t}}\biggl[\exp\biggl\{-\frac{(y-x)^2}{2t}\biggr\}\vspace*{1pt}
\cr
\displaystyle \hspace*{66.5pt}{}-\exp\biggl\{-\frac{(y+x)^2}{2t}\biggr\}
\biggr]\delta_0(dl) \,dy,
\qquad\mbox{if } y\ge0,l \ge0,}
\end{eqnarray*}
whereas if $x \leq0$, then
\begin{eqnarray*}
&&
P_x\bigl(B_t^{(\alpha)}\in dy, \ell_t^{(\alpha)}\in d\ell\bigr) \\
&&\qquad=\cases{
\displaystyle \frac{2\alpha(l+y-x)}{\sqrt{2\pi t^3}}\exp\biggl\{-\frac{(l+y-x)^2}{2t}\biggr\}
\,dy\,dl,\vspace*{5pt}\cr
\hspace*{186pt}\qquad\qquad\mbox{if } y\ge0,l \ge0 ,\vspace*{5pt}\cr
\displaystyle \frac{2(\alpha-1)(l-y-x)}{\sqrt{2\pi t^3}}\exp\biggl\{-\frac{(l-y-x)^2}{2t}\biggr\}\,dy\,dl\vspace*{5pt}\cr
\qquad{}+ \displaystyle \frac{1}{\sqrt{2\pi
t}}\biggl[\exp\biggl\{-\frac{(y-x)^2}{2t}\biggr\}\vspace*{5pt}\cr
\displaystyle \hspace*{66.5pt}{}-\exp\biggl\{-\frac{(y+x)^2}{2t}\biggr\}
\biggr]\delta_0(dl) \,dy,
\qquad\mbox{if } y\leq0,l \ge0.}
\end{eqnarray*}
\end{cor}

\section*{Acknowledgments}

The authors are grateful to Pierre Etor\'e and Miguel Martinez
for sharing their preprint \citet{etoiremartinez}, who pointed out
these mistakes.
In particular, the formula with $(2\alpha-1)v$ in place
of $|(2\alpha-1)v|$ is given in
Proposition 4.1 of their paper.



%
\printaddresses

\end{document}